\pdfoutput=1 %arXiv admin added this line
\documentclass{article}

\usepackage{amsmath}
\usepackage{amssymb}
\usepackage{graphicx}

\evensidemargin \oddsidemargin
\newtheorem{theorem}{Theorem}[section]
\newtheorem{lemma}[theorem]{Lemma}

\def\F{\mathcal F}
\def\G{\mathcal G}
\def\H{\mathcal H}

\def\O{\mathcal O}
\def\R{\mathcal R}

\begin{document}

\begin{center}
{\large\bf Numerical Smoothness and Error Analysis for RKDG \\
on the Scalar Nonlinear Conservation Laws}

\vskip.20in

Tong Sun   and   David Rumsey   %%$^{1}$
%\\[2mm]

%%$^{1}$
Department of Mathematics and Statistics\\
Bowling Green State University \\
Bowling Green, OH 43403

\end{center}

\date{today}

%{\footnotesize

%\begin{article}
%\keywords{Long-time error estimation, stability-smoothing indicator, parabolic PDE, %numerical solution, subdomain, partially implicit scheme.}

%{\bf AMS subject class.} Primary: 65M15  Secondary: 65L20

%\authorrunninghead{SUN}
%\titlerunninghead{Consistency + Numerical Smoothing $\Rightarrow$ Convergence}

\setcounter{page}{1}   %% This command is optional.
                       %% May set page number only for first page in
                       %% issue, if desired.

%% &lt;&lt;== End of commands to be entered at Wiley

%%  Authors, start here ==&gt;&gt;

%\title{Consistency + Numerical Smoothing $\Rightarrow$ Convergence \\
%   \hspace{0.6in} {\it An Alternative of the Lax-Richtmyer Theorem }}

%\author{Tong Sun}
%\affil{Department of Mathematical and Statistics\\
%Bowling Green State University\\
%Bowling Green, Ohio, USA}

\abstract{ The new concept of numerical smoothness is applied to the
RKDG (Runge-Kutta/Discontinuous Galerkin) methods for scalar
nonlinear conservations laws. The main result is an {\it a
posteriori} error estimate for the RKDG methods of arbitrary order
in space and time, with optimal convergence rate. In this paper, the
case of smooth solutions is the focus point. However, the error
analysis framework is prepared to deal with discontinuous solutions
in the future. }

\section{Introduction}

The Runge-Kutta/Discontinuous Galerkin (RKDG) methods are among the
most popular modern numerical methods for nonlinear conservation
laws. Due to the complexity of the schemes and the nonlinearity of
the problems, error analysis theory for the RKDG methods is not yet
satisfactorily completed. A brief summary of the currently available
error analysis results can be found in the recent papers by Q. Zhang
and C.-W. Shu, \cite{qz2004} and \cite{qz2009}. Here in this paper,
a new error analysis is given based on the innovative concept of
numerical smoothness. The main result of this paper is a practical
{\it a posteriori} error estimate of optimal order , which depends
on a set of computed smoothness indicators.

In most {\it a priori} error analysis of time-dependent problems,
local error is referred to original PDE solutions to take advantage
of their smoothness, consequently error propagation has to be
referred to numerical schemes (e.g. \cite{qz2009}). In {\it a
posteriori} error analysis, local error can only be referred to
numerical solutions, consequently error propagation is referred to
PDEs. Since numerical solutions of PDEs are typically not smooth
functions (discrete point values in finite difference schemes,
piecewise polynomials in finite elements, etc.), when local error is
referred to numerical solutions, it is often given as residuals,
such as in the well-known duality method \cite{cockB}
\cite{johnson}. In our {\it a posteriori} error analysis of the
scalar conservation laws and RKDG, we use the $L_1$-contraction
between the PDEs' entropy solutions for error propagation analysis,
and rely on numerical smoothness instead of residuals to estimate
local error.

The idea of using numerical smoothness in the error analysis of
nonlinear hyperbolic conservation laws is a migration of the idea of
using numerical smoothing in the error analysis of nonlinear
parabolic equations solved with complex schemes
\cite{sun}\cite{sunF}. For nonlinear equations solved with complex
schemes, we base the concept of numerical smoothness on a set of
efficiently computable smoothness indicators. When the indicators
remain bounded during actual computation, we consider the numerical
solution as being numerically smooth. Due to the equations'
nonlinearity and the schemes' complexity, we usually cannot give an
{\it a priori} proof the boundedness of our smoothness indicators.
However, we can always compute the indicators along a numerical
solution. We define the indicators at $t_n$ in such a way that we
can prove a local error estimate of optimal rate for the time step
$[t_n,t_{n+1}]$, where the indicators play the role of high order
derivatives as in most {\it a priori} error estimates. It will be
shown that numerical smoothness indicators deliver much more
abundant information then residuals. Consequently, we can get better
error estimates and more useful information toward adaptive
algorithms.

Numerical solutions are not always numerically smooth. The usual
measures taken for the purpose of achieving numerical stability are
actually also what is needed to achieve numerical smoothness,
because one of the main ingredients of these measures is numerical
diffusion. Smoothness indicators serve two purposes: (1) to watch
the smoothing and/or smoothness maintenance performance of the
scheme; and (2) to provide smoothness information for local error
estimates. For the RKDG schemes, we use the Godunov upwind flux, the
TVD-RK schemes and a strengthened CFL condition. After taking all
these measures, it is extremely hard to prove the boundedness of our
smoothness indicators. However, the importance of the whole idea
resides on the fact that we can use the computed smoothness
indicators to circumvent the difficult proof of numerical
smoothness, and move on to prove sharp error estimates. For complex
nonlinear problems, this kind of circumvention is probable the only
way to achieve practical error estimates.

While doing the proofs, we realized another advantage of the
numerical smoothness approach. Since we are working on DG finite
element solutions and the entropy solutions just evolving away from
those finite element solutions, we have easy access to the
$L_\infty$-norm estimates, which in turn gives us $L_1$ and $L_2$
estimates. For error propagation, $L_1$ contraction is the best
tool. For the finite element formulations, $L_2$-norm is natural.
For the nonlinearity, $L_\infty$-estimates are crucial. Having
access to all three, we are able to do the error analysis of the
finite element methods of the nonlinear problems, where the global
error estimates do not have an exponentially growing factor.

In the solutions of nonlinear conservation laws, there are shocks
and contact discontinuities. In the discontinuous Galerkin finite
elements, there are the technical discontinuities of the piecewise
polynomials. The analysis of this paper is limited to the case of
smooth PDE solutions; henceforth, only the technical discontinuities
are treated. While we do not assume anything directly on the
smoothness of the PDE's solution, we only consider the case that all
the components of our smoothness indicators are well-bounded. In
fact, the boundedness of the smoothness indicators indicates that a
smooth PDE solution is being approximated.  Our smoothness
indicators are capable of detecting shocks and contact
discontinuities (including high order discontinuities)
\cite{rumsey}. Our $L_1$-contraction error propagation analysis
remains valid in dealing with shocks and other discontinuous
solutions. However, we restrict ourselves to the case of smooth
solutions in this paper. Many discontinuous solutions of
conservation laws are piecewise smooth. The work of this paper can
also be considered as analyzing error in a smooth piece of a
discontinuous solution. Clearly, we need to understand how to obtain
optimal error estimates on smooth pieces of solutions, before we
focus on the error analysis at shocks and contact discontinuities.
In this sense, this paper is the first step of the project of
analyzing the error of RKDG methods with numerical smoothness.

In this paper, our goal is to show the new error analysis ideas and
the nature of the results. In order to focus on the framework, we do
not trace all the constants involved in the error estimates.
Instead, we show how they should be computed with enough details to
reveal their dependency, computability and boundedness. A separate
technical report will be prepared to show the fine details. Since
some of the constants depend on the flux function $f$, certain
details are better shown with numerical examples. No generic
constant will appear in this paper.

The nonlinear conservation law problems and the RKDG schemes are
well-known. For a survey article, see the lecture notes \cite{ShuD}
by C.-W. Shu.  Consider the one-dimensional nonlinear conservation
law
\begin{equation}
u_t + f(u)_x = 0 \label{law}
\end{equation}
in a bounded interval $\Omega = [a,b]$.  In order to focus on the new ideas and the
new tools of the proof, we stay with the simple case of west wind ( $f'(u) > 0$ ) .
Let the initial condition be
\begin{equation} \label{initialU}
u(0,x) = u_I(x)
\end{equation}
and the upwind boundary condition be
\begin{equation} \label{boundU}
u(t,a) = u_L(t).
\end{equation}
Assume that the flux function $f(u)$ is sufficiently smooth and the
initial and boundary conditions are smooth and consistent to
guarantee that the entropy solution $u(t,x)$ is smooth near $x=a$
for all $t>0$.

Partition $\Omega$ with $a=x_{-1/2} < x_{1/2} < \cdots < x_{m-1/2} =
b$. Let $h = x_{j+1/2} - x_{j-1/2}$ be the same for all cells
$\Omega_j = [ x_{j-1/2} , x_{j+1/2}]$. To solve the problem with the
discontinuous Galerkin method, we take the standard discontinuous
piecewise polynomials space $V_h$. When the degree of a local
polynomial is up to $p$, $V_h = \{ v \in L^2(\Omega) : v|_{\Omega_j}
\in \Pi_p \}$, where $\Pi_p$ is the set of all the polynomials of
degree less than or equal to $p$. In each cell $\Omega_j$, a
semi-discrete solution $u_h $ satisfies
\begin{equation}
(u_{h,t}, v)_{\Omega_j} = (f(u_h), v_x)_{\Omega_j}
                        + f(u_h(x_{j-1/2}^-)) v(x_{j-1/2}^+)
                        - f(u_h(x_{j+1/2}^-)) v(x_{j+1/2}^-). \label{semi}
\end{equation}
Here the Godunov flux is employed under the west wind assumption
(for simplicity). At the upwind boundary $x_{-1/2}=a$, we set
\begin{equation} \label{boundUh}
u_h(t,x_{-1/2}) = u(t,a) = u_L(t).
\end{equation}
At the initial time $t=0$, $u_h(0,x)$ is taken to be the
$L_2$-projection of $u_I(x)$.

For temporal discretization, we take a standard TVD-RK scheme of
order $k$ \cite{ShuD}. For example, when $k=3$, in the time step
$[t_n, t_{n+1}]$, we compute the fully-discrete solution $u^c_{n+1}
\in V_h$ from $u^c_n$ by the following scheme. With the notation
$$
\H_j(u,v) = (f(u),v_x)_{\Omega_j} + f(u(x_{j-1/2}^-)) v(x_{j-1/2}^+)
- f(u(x_{j+1/2}^-)) v(x_{j+1/2}^-) ,
$$
the scheme is that, for $\tau_n = t_{n+1} - t_n$ and any $v \in
V_h$,
\begin{equation} \label{stage1}
(u^{c,1}_n,v)_{\Omega_j}  = (u^{c}_n,v)_{\Omega_j}   + \tau_n \H_j(u^c_n, v),
\end{equation}
\begin{equation} \label{stage2}
(u^{c,2}_n,v)_{\Omega_j}  = \frac{3}{4} (u^{c}_n,v)_{\Omega_j}
                                        +  \frac{1}{4} (u^{c,1}_n,v)_{\Omega_j}  + \frac{\tau_n}{4} \H_j(u^{c,1}_n, v),
\end{equation}
\begin{equation} \label{stage3}
(u^{c}_{n+1},v)_{\Omega_j}  = \frac{1}{3} (u^{c}_n,v)_{\Omega_j}
                                            +  \frac{2}{3} (u^{c,2}_n,v)_{\Omega_j}  + \frac{2 \tau_n}{3} \H_j(u^{c,2}_n, v).
\end{equation}
At $t=0$, we make $u^c_0 = u_h(0,x)$. The upwind boundary values of
$u^{c,1}_n$, $u^{c,2}_n$, and $u^c_{n+1}$ are taken according to
$u_L(t)$.

\section{Error propagation and numerical smoothness}

The PDE's entropy solution satisfying the original initial condition
is denoted by $u(t,x)$.   Throughout this paper, we only consider
one numerical solution, namely the computed numerical solution,
which is denoted by $u^c_n$ as above. In order to present the local
error analysis in a time step $[t_n,t_{n+1}]$, we use the
semi-discrete solution that passes $(t_n,u^c_n)$.  For the briefness
of notations, we use $u^{h}(t,x)$ for this temporally piecewise
semi-discrete solution, which has a new initial value in each time
step. $u^h$ takes the upwind boundary value given in
(\ref{boundUh}). Since we do not simultaneously work on the local
error analysis of two different time steps, the notation
$u^{h}(t,x)$ should not cause ambiguity. In each time step, we also
need the PDE's entropy solution which passes $(t_n,u^c_n)$. We
denote this entropy solution by $\tilde{u}(t,x)$. $\tilde{u}(t,x)$
satisfies the upwind boundary condition (\ref{boundU}). Of course,
$\tilde{u}$ is also defined piecewise in time. The following error
splitting diagram may help the reader in remembering the notations
for these solutions.

\setlength{\unitlength}{1.0mm}
\begin{picture}(70,30)(-10,10)
    \linethickness{0.6pt}
    \qbezier(30,8)(50,14)(70,26)
    \linethickness{0.5pt}
    \qbezier(30,0)(50, 6)(70,18)
    \qbezier(30,0)(50, 6)(70,15)
    \qbezier(30,0)(50, 6)(70,12)
    \put(25,-2){\makebox(0,0)[b]{$u^c_n$}}
    \put(24, 6){\makebox(0,0)[b]{$u(t_n)$}}
    \put(78, 9){\makebox(0,0)[b]{$u^c_{n+1}$}}
    \put(80.4,13.5){\makebox(0,0)[b]{$u^{h}(t_{n+1})$}}
    \put(80,19){\makebox(0,0)[b]{$\tilde{u}(t_{n+1})$}}
    \put(80,25){\makebox(0,0)[b]{$u(t_{n+1})$}}
    \put(5, -10){\vector(1,0){85}}
    \put(10,-11){\vector(0,1){45}}
    \put(29.8,-10){\line(0,1){18}}
    \put(70,-10){\line(0,1){36}}
    \put(30,-14){\makebox(0,0)[b]{$t_n$}}
    \put(70,-14){\makebox(0,0)[b]{$t_{n+1}$}}
\end{picture}

\vspace{1.2in}

In the diagram and also in the rest of the paper, sometimes we hide
one of the two independent variables in the notation of a solution
to make the expressions shorter.

The error analysis of this paper is based on the error splitting in
the diagram. In order to estimate the global error
$u(t_{n+1})-u^c_{n+1}$ at time $t_{n+1}$, we split it into three
parts as shown in the diagram.
\begin{eqnarray}
\|u(t_{n+1})-u^c_{n+1}\|
&\le& \|u(t_{n+1}) - \tilde{u}(t_{n+1})\|  \nonumber \\
&+& \|\tilde{u}(t_{n+1}) - u^{h}(t_{n+1})\| \nonumber \\
&+& \|u^{h}(t_{n+1}) - u^c_{n+1}\| \label{split}
\end{eqnarray}

The first part $u(t_{n+1}) - \tilde{u}(t_{n+1})$ is the propagation
of the global error $u(t_n) - u^c_n$ by the PDE. Due to the
$L_1$-contraction property of the scalar conservation laws, since
$u$ and $\tilde{u}$ satisfy the same upwind boundary condition
(\ref{boundU}), we have
\begin{equation} \label{contraction}
\|u(t_{n+1}) - \tilde{u}(t_{n+1})\|_{L_1(\Omega)} \le  \|u(t_n) -
u^c_n\|_{L_1(\Omega)}.
\end{equation}

The second part of the split error is the local spatial
discretization error $\tilde{u}(t_{n+1}) - u^{h}(t_{n+1})$. Since
$u^c_n$ lives in a discontinuous finite element space, $\tilde{u}$
is certainly not smooth in the classical sense. In fact, the
discontinuity of $\tilde{u}(t_n) = u^c_n$ at $x_{j-1/2}$ will travel
into the cell $\Omega_j$. Thus, at any time $t \in (t_n,t_{n+1}]$,
there is either a shock, a contact discontinuity, or a rarefaction
wave of $\tilde{u}$ in $\Omega_j$. However, if the solution $u(t,x)$
is smooth around $\Omega_j$, we intuitively know that the
discontinuity of $\tilde{u}$ is only technical and it must be very
tiny, and $\tilde{u}$ must be smooth away from its discontinuities.
We will quantitatively substantiate the intuition in the definition
of the spatial smoothness indicator. Then we will use the indicator
to estimate the spatial local error.

The third part of the split error is the local temporal
discretization error $u^{h}(t_{n+1}) - u^c_{n+1}$. To estimate
this part of the error, we will need the temporal smoothness of
$u^{h}$ for $t \in [t_n,t_{n+1}]$. Since $u^{h}$ is an ODE
solution with initial value $u^c_n$, we can also establish the
needed smoothness.

Here it is to be noticed that the analysis relies on the smoothness
properties of $\tilde{u}$ and $u^{h}$. Since both of them have
$u^c_n$ as their initial value, the smoothness level of them depends
on how $u^c_n$ has been computed. In other words, some kind of
numerical smoothing or smoothness-maintenance should have been built
in the scheme. In this paper, we do not intend to prove such
smoothing or smoothness-maintenance ability for the RKDG methods.
Fortunately, the computed smoothness indicators in our numerical
experiments show that the RKDG methods do have the desired ability
to keep a numerical solution ``smooth" (when/where the solution
should be smooth). We only prove error estimates by using the
smoothness indicators.

\section{The smoothness indicators}
In order to rigorously and quantitatively define the concept of
numerical smoothness for the RKDG method, we define the following
spatial and temporal smoothness indicators for each time step
$[t_n,t_{n+1}]$.
\begin{itemize}
\item Spatial smoothness indicator: $S_n^{p} = S_{p}(u^c_n),$
\item Temporal smoothness indicator: $T_n^{k} = T_{k}(u^c_n).$
\end{itemize}
Here $S$ stands for space, $T$ stands for time, $p$ is the degree of
the polynomials in each cell, and $k$ is the order of the
Runge-Kutta scheme.

\subsection{Definition of $T^k_n$}

The temporal smoothness indicator $T_n^{k}$ consists of the temporal
derivatives of $u^{h}$ at $t=t_n$. Namely,
$$ T_n^{k} = (u^c_n, u^{h}_{t}(t_n), u^{h}_{tt}(t_n),
\cdots, \frac{\partial^{k+1} u^{h}}{\partial t^{k+1}}(t_n)).$$ The
first derivative $u^{h}_{t}(t_n)$ is computed as in the
implementation of the forward Euler scheme
$$
(u^{h}_{t}(t_n),v)_{\Omega_j} = \H_j(u^c_n, v).
$$
The formula for computing the second derivative can be obtained by
taking derivatives with respect to time on both sides of the
semi-discrete scheme (\ref{semi}):
\begin{eqnarray*}
(u^h_{tt}, v)_{\Omega_j} &=& (f'(u^h) u^h_{t}, v_x)_{\Omega_j} \\
        &+& f'(u^h(x_{j-1/2}^-)) u^h_{t}(x_{j-1/2}^-) v(x_{j-1/2}^+) \\
        &-& f'(u^h(x_{j+1/2}^-)) u^h_{t}(x_{j+1/2}^-)
        v(x_{j+1/2}^-).
\end{eqnarray*}
To compute $u^{h}_{tt} (t_n)$ with this formula,  on the right hand
side, we replace $u^h$ by $u^c_n$ and  $u^h_{t}$ by the computed
first derivative $u^{h}_{t}(t_n)$. The high order derivatives can be
computed similarly.

The ability of the indicator $T_n^{k}$ to reveal numerical
solutions' smoothness, discontinuities, and possible numerical
``instability" phenomena has been reported in \cite{rumsey}. Since
$T_n^{k}$ contains the initial temporal derivatives of $u^{h}$, it
can be used for the temporal local error estimation without any
transformation.

\subsection{Definition of $S^p_n$}

The spatial smoothness indicator $S^p_n$ contains not only the
spatial derivatives of $u^c_n$ within each cell, but also the jumps
of the derivatives across the cell boundaries. Namely, for the cell
$\Omega_j = [x_{j-1/2} ,x_{j+1/2}]$,
$$ S_{n,j}^{p} = (M_{n,j}^0, M_{n,j}^1, \cdots, M_{n,j}^p, D_{n,j}^0, D_{n,j}^1, \cdots, D_{n,j}^p), $$
where
$$ M_{n,j}^l = \frac{\partial^l}{\partial x^l} u^c_n(x_{j-1/2}^+), \qquad
L_{n,j}^l = \frac{\partial^l}{\partial x^l} u^c_n(x_{j-1/2}^-),$$ and
$$ J_{n,j}^l = M_{n,j}^l - L_{n,j}^l = D_{n,j}^l {h^{p+1+\mu -l(1+\alpha)}} , $$
for some properly determined constants $\alpha \in [0,1)$ and $\mu
\in [0,1]$.

When $j=0$, $L_{n,0}^l$ needs to be defined separately. To this end,
we use the upwind boundary function $u(t,x_{-1/2}) = u_L(t)$. By
calculation from the conservation law (\ref{law}), we must have
$$
L^0_{n,0} = u(t_n,x_{-1/2}) = u_L(t_n),
$$
$$
L^1_{n,0} = u_x(t_n,x_{-1/2}) =
-\frac{\frac{d}{dt}u_L(t_n)}{f'(u_L(t_n))},
$$
$$
L^2_{n,0} = u_{xx}(t_n,x_{-1/2}) = -\frac{2 f''(u_L(t_n))
[\frac{d}{dt}u_L(t_n)]^2 - f'(u_L(t_n))
\frac{d^2}{dt^2}u_L(t_n)}{[f'(u_L(t_n))]^3},
$$
and so on. For later use, we extend $u^c_n$ to $\Omega_{-1} =
[x_{-1/2}-h,x_{-1/2}]$ by $u^c_n = u(t_n,x)$, where $u(t_n,x)$ is
obtained by a short time tracing back from $u(t,a)=u_L(t)$.  Under a
proper smoothness assumption on $u_L(t)$, such tracing back is well
defined. By Taylor expansion,
$$
u^c_n(x) = u(t_n,x) = L_{n,0}^0 + L_{n,0}^1 (x-x_{-1/2}) + \cdots +
L_{n,0}^p (x-x_{-1/2})^p/{p!} + R_{n,0}(x), \qquad x\in\Omega_{-1}.
$$
The residual $R_{n,0}(x) = \O((x-x_{-1/2})^{p+1})$ is of higher
order. Given a smooth boundary function $u_L$, one can determine a
constant $\bar{D}$, such that
$$
| R_{n,0}(x) | \le \bar{D} |x-x_{-1/2}|^{p+1} /(p+1)!
$$
In other cells ($j>0$), let $R_{n,j}(x) = 0$.  The expansion part of
$u^c_n$ is not computable, it only lives in the proof. The expansion
is defined in this way to be consistent with the boundary condition
satisfied by $\tilde{u}$.

It is obvious that the values of $M^l_{n,j}$ and $L^l_{n,j}$ should
be of $\O(1)$, unless there is a shock or contact discontinuity
somewhere around $\Omega_j$. It is also easy to guess that the jumps
$J_{n,j}^l$ should be small, otherwise the numerical solution may
have lost too much smoothness around the cell boundary. How small
should the jumps $J_{n,j}^l$ be? Both our error analysis and
numerical experiments suggest that $D_{n,j}^l =
J_{n,j}^l/h^{p+1+\mu-l(1+\alpha)}$ should be at most of $\O(1)$,
unless there is a shock or high order discontinuity within or near
the cell. This is the reason for having $D_{n,j}^l$ instead of
$J_{n,j}^l$ serving as a part of the smoothness indicator.

How is $\alpha$ determined? It is well known that, with high
degree DG elements, the time step size $\tau_n$ should satisfy a
strengthened CFL condition of the form
$$ \tau_n < \gamma h^{1+\alpha}.$$
In \cite{qz2004}, for example, $\alpha = 1/3$. For the definition of
$D^l_{n,j}, l=0,\cdots,p$, we need $\alpha = \mu / p$. In fact,
since $J^p_{n,j}$ is the jump of the piecewise constant function
$\frac{\partial^p}{\partial x ^p} u^c_n$, to match the total
variation of $\frac{\partial^p}{\partial x ^p} u$, the average of
the jumps $J^p_{n,j}$ must be of $\O(h)$. That is, we need
$p+1+\mu-p(1+\alpha) = 1$ in the definition of $D^p_{n,j}$, or
equivalently $\mu = p \alpha$, to have $D^p_{n,j} = \O(1)$.

The amount of work for computing $S^p_n$ is proper. In fact, $S^p_n$
contains the minimal amount of smoothness information for us to
estimate the optimal approximation error of the piecewise
polynomials of degree $p$ to the PDE solution. The amount of work
for computing $T^k_n$ is also proper, for the same reason. It might
be possible to estimate $T^k_n$ from $S^p_n$ (if $k$ and $p$ are
related in certain way), but a directly computed $T^k_n$ should
sharpen the temporal local error estimate.

We consider the numerical solution as a good approximation of a
smooth PDE solution if and only if $S^p_n$ is reasonably bounded. It
will be a future issue to study how to classify the numerical
solution in case $S^p_n$ is not reasonably bounded. We will have to
distinguish different patterns of $S^p_n$. What indicates a
well-caught shock, or well-approximated transition to a shock? What
indicates a well-approximated high order contact discontinuity? What
indicates numerical ``instability"? How to adaptively deal with each
of these cases? The temporal smoothness indicator $T^k_n$ should
also be studied for the same issues.

\section{The main error estimates}

\begin{theorem} \label{main}
Let $u(t,x)$ be the entropy solution of the nonlinear conservation
law (\ref{law}) satisfying the initial condition (\ref{initialU})
and upwind boundary condition (\ref{boundU}). Let $u^c_n$ be the
numerical solution computed by a TVD-RK-DG scheme with piecewise
polynomials of degree $p$ and the TVD-RK scheme of order $k$, on the
partition of $\Omega$ described in Section 1. Assume that $u$ and
$u^c_n$ are bounded by a constant $U$ in $[0,T] \times \Omega$. Let
$\beta = \max_{|w| \le U} f'(w)$. Assume that the time step size
$\tau_n$ for each step satisfies the standard CFL condition $\beta
\tau_n \le h$ and the strengthened CFL condition $\tau_n \le \gamma
h^{1+\alpha}$, for a constant $\mu \in [0,1]$, a positive constant
$\gamma$, and $\alpha = \mu/p$.

If there is a positive real number $M$, such that, for all $t_n \le
T$, all the components of $S^p_n$ and $T^k_n$ are bounded by $M$,
then the spatial and temporal local error in $[t_n,t_{n+1}]$ satisfy
\begin{equation}
\|\tilde{u}(t_{n+1}) - u^{h}(t_{n+1})\|_{L_1(\Omega)} \le \tau_n
h^{p+\mu}  \F (S_n^{p}), \label{spaceError}
\end{equation}
\begin{equation}
\|u^{h}(t_{n+1}) - u^c_{n+1}\|_{L_1(\Omega)} \le \tau_n^{k+1}
\G(T_n^{k},S_n^{p}), \label{timeError}
\end{equation}
where $\F(S^p_n)$ and $\G(T^k_n,S_n^{p})$ are computable functions
of the indicators. As a consequence of the error splitting
(\ref{split}), the $L_1$-contraction property (\ref{contraction}),
and the local error estimates (\ref{spaceError}) and
(\ref{timeError}),
$$
\|u(t_{n+1}) - u^c_{n+1}\|_{L_1(\Omega)} \le \|u(t_{n}) -
u^c_{n}\|_{L_1(\Omega)} + \tau_n [ h^{p+\mu} \F(S_n^{p}) +  \tau_n^k
\G(T_n^{k},S_n^{p}) ].
$$
Finally, at the end of the computation ($t_N = T$),
\begin{equation} \label{global}
\|u(T) - u^c_N\|_{L_1(\Omega)} \le \|u(0) - u^c_0\|_{L_1(\Omega)} +
\Sigma_{n=1}^N \tau_n [ h^{p+\mu} \F(S_n^{p}) +  \tau_n^k
\G(T_n^{k},S_n^{p}) ].
\end{equation}
\end{theorem}

\noindent {\bf Proof.} It suffices to prove (\ref{spaceError}) and
(\ref{timeError}).  The next two subsections will carry the proofs
of (\ref{spaceError}) and
(\ref{timeError}) respectively. \# \\

{\bf Remark.} In the literature, $\mu=1$ is considered to be the
optimal convergence rate. We keep $\mu$ as a parameter to cover
those possible non-optimal cases. However, when the initial solution
is smooth, we always have $\mu=1$. For $p \ge 3$, $\alpha = \mu/p$
is not too restrictive. There is no restriction on $\gamma$ in the
proof, although $\gamma$ will appear in the function $\F(S^p_n)$.
The actual restriction on $\tau_n$ is in real computation. If
$\tau_n$ is too large, the RKDG scheme fails on numerical smoothness
maintenance. See the numerical experiments in Section 5.

\subsection{Estimating $\tilde{u}(t_{n+1}) - u^{h}(t_{n+1})$, proof
of (\ref{spaceError}) }

We begin with introducing an auxiliary piecewise PDE solution
$u^{e}$. First define a local strong solution $u^{e}_j$ of the
conservation law. The initial values of $u^e_j$ are given on the
line segment $\{t_n\} \times (\Omega_{j-1} \cup \Omega_j)$ by
$$
u^{e}_j(t_n,x) = M_{n,j}^0 + M_{n,j}^1 (x-x_{j-1/2}) + \cdots +
M_{n,j}^p (x-x_{j-1/2})^p/{p!} \qquad  j=0,1,\cdots,m-1.
$$
It is easy to see that, in $\Omega_j$, $u^{e}_j(t_n) = u^c_n$; in
$\Omega_{j-1}$,
$$
u^{e}_j(t_n) =u^c_n + J_{n,j}^0 + J_{n,j}^1 (x-x_{j-1/2}) + \cdots +
J_{n,j}^p (x-x_{j-1/2})^p/{p!} - R_{n,j}(x).
$$
As a strong solution of the Cauche problem of the original
conservation law (\ref{law}),  $u^{e}_j$ certainly exists in the
region $\R_{n,j} = \{(t,\tilde{x})| t \in [t_n,t_{n+1}], x \in
\Omega_{j-1} \cup \Omega_j, \tilde{x} \le x_{j+1/2}, \tilde{x} = x +
f'(u^{e}_j(t_n,x))(t-t_n) \}$. This is the trapozoidal region
covered by the characteristic lines originating from $\Omega_{j-1}
\cup \Omega_j$. When $\tau_n$ satisfies the standard CFL condition
$\beta \tau_n \le h$, it is easy to see that $[t_n, t_{n+1}] \times
\Omega_j \subset \R_{n,j} \subset [t_n, t_{n+1}] \times
(\Omega_{j-1} \cup \Omega_j)$.

At the upwind boundary, let $u^{e}_{-1} = u$ for $x \in \Omega_{-1}
= [x_{-1/2}-h,x_{-1/2}]$. Due to the smoothness of $u_L(t)$, one can
determine the value of $u$ in $\Omega_{-1}$ by tracing back (but not
computable).  Now, we are ready to define the local piecewise PDE
solution by
$$
u^{e}(t,x) = u^{e}_{j}(t,x), \qquad (t,x) \in [t_n,t_{n+1}] \times \Omega_j, \, \,\,
j=-1,0,1,\cdots,m-1.
$$
Since $u^{e}_j(t_n)$ is a polynomial in $\Omega_{j-1} \cup
\Omega_j$, $u^{e}$  is smooth in $[t_n,t_{n+1}] \times \Omega_j$ for
sufficiently small $\tau_n$. To reveal more details on the
smoothness of $u^e$, we have the following Lemma.

\begin{lemma} \label{lemma1}
 There are constants $N^l_{n,j}$  ($l=0,1,\cdots, p+1,  j=0,1,\cdots,m-1$), which depend on the flux
 function $f$ and can be computed from $M^0_{n,j}, M^1_{n,j}, \cdots, M^p_{n,j}$, such
 that,
$$
\| \frac{\partial^l}{\partial x^l} u^{e}_{j}(t,x)
\|_{L_{\infty}(\R_{n,j})} \le N^l_{n,j}, \qquad l = 0,1,\cdots, p,
\,\,\, j=0,1,\cdots,m-1.
$$
Moreover,
$$
\| \frac{\partial^{p+1}}{\partial x^{p+1}} u^{e}_{j}(t_n+\tau,x)
\|_{L_{\infty}(\Omega_j)} \le \tau N^{p+1}_{n,j}, \qquad
j=0,1,\cdots,m-1 .
$$
\end{lemma}
\noindent {\bf Proof.} For the simplicity of notations, in the
proofs of this  and the next Lemma, we denote the solution
$u^{e}_{j}$ by $w$, and $\frac{\partial^l w}{\partial x^l}$ by
$w^{(l)}$. From
$$
w_t + f(w)_x = 0,
$$
we get
$$
w^{(1)}_t + f'(w) w^{(1)}_x + f''(w) ( w^{(1)} )^2 = 0,
$$
$$
w^{(2)}_t + f'(w) w^{(2)}_x + 3 f''(w) w^{(2)} w^{(1)} + f'''(w) (w^{(1)})^3 = 0 ,
$$
and so on. The boundedness of $w$ is obvious.

Along each characteristic line, $w^{(1)}$ will not have a blow-up in
a short  time $\tau \in [0,\tau_n]$, where $\tau = t-t_n$. The
initial value of $w^{(1)}$ is $\frac{\partial}{\partial x}
u^{e}_{j}(t_n) = M_{n,j}^1 + M_{n,j}^2 (x-x_{j-1/2}) + \cdots +
M_{n,j}^p (x-x_{j-1/2})^{p-1}/{(p-1)!}$. $f''(w)$ can also be
computed from $u^{e}_{j}(t_n)$. Therefore, $w^{(1)}$ can be
estimated by the values of $M^0_{n,j}, M^1_{n,j}, \cdots, M^p_{n,j}$
in the $L_\infty$ norm. Hence $N^1_{n,j}$ can be computed.

As for the higher order derivatives, the ODE for each $w^{(l)}$
along each characteristic line is linear in $w^{(l)}$ and depends on
the lower order derivatives $w, w^{(1)}, \cdots, w^{(l-1)}$. The
initial value of $w^{(l)}$ only depends on $u^{e}_j(t_n)$,
therefore the bound $N^l_{n,j}$ of $w^{(l)}$ can also be estimated
from $M^0_{n,j}, M^1_{n,j}, \cdots, M^p_{n,j}$, as a result of
mathematical induction.

$w^{(p+1)}$ has a special property: its initial value is the 0
function (the $p$+1-st derivative of $u^{e}_j(t_n)$). Integrating
the ODEs about $w^{(p+1)}$ along each characteristic line, we
realize that $w^{(p+1)}$ is proportional to $\tau$, while the
coefficient $N^{p+1}_{n,j}$  can be computed from
$M^0_{n,j}, M^1_{n,j}, \cdots, M^p_{n,j}$.    \#  \\

{\bf Remarks.} Here we make two remarks on the results of Lemma
\ref{lemma1}. (1) Because each ODE is integrated along a
characteristic line for a short time (not longer than a time step),
one can simply use $M^l_{n,j}$ as a practical estimate for
$w^{(l)}$. In other words, $N^l_{n,j}$ is actually very close to
$M^l_{n,j}$. However, the most useful $N^{p+1}_{n,j}$ has to be
computed through solving the differential inequalities. (2) The
factor $\tau$ in the estimate of $w^{(p+1)}$ is crucial for the
error analysis later on. It means that, because $u^{e}_{j}$ is
evolving out of a polynomial of degree $p$, when $u^{e}_{j}$ is
approximated by a polynomial of degree $p$, the error is
proportional to $\tau$. The idea of picking up this $\tau$ for the
local spatial error (in one way or another) comes from reading \cite{qz2009}. \\

As it appears in most error analysis of finite element methods, we
also need an $L_2$-projection of a smooth solution. To this end, we
consider the cell by cell $L_2$-projection of $u^{e}$ into $V_h$.
Denote this projection by $u^{p} =u^{p}(t,x) \in V_h$ (p
stands for projection here), it is given by
$$
(u^{p},v)_{\Omega_j} = (u^{e}, v)_{\Omega_j}, \qquad \forall v
\in V_h.
$$
By the Bramble-Hilbert Lemma, the scaling argument, and Lemma \ref{lemma1}, the following
estimates are obvious:
\begin{lemma} \label{lemma2}
For sufficiently small $\tau = t-t_n$,
$$
\|u^e - u^p\|_{L_1(\Omega_j)} \le C_1 h^{p+1} \|w^{(p+1)}\|_{L_1(\Omega_j)}  \le C_1 h^{p+1} \tau h N^{p+1}_{n,j},
$$
$$
\|u^e - u^p\|_{L_2(\Omega_j)}  \le C_2 h^{p+1} \|w^{(p+1)}\|_{L_2(\Omega_j)} \le C_2 h^{p+1} \tau\sqrt{h}  N^{p+1}_{n,j},
$$
and
$$
\|u^e-u^p\|_{L_{\infty}(\Omega_j)}  \le C_3 h^{p+1} \|w^{(p+1)}\|_{L_{\infty}(\Omega_j)} \le C_3 h^{p+1} \tau N^{p+1}_{n,j}.
$$
Here $C_1, C_2$ and $C_3$ are the projection error constants in the
reference  cell. Consequently, in the whole domain $\Omega$, let
$N^{p+1}_n = \max_j N^{p+1}_{n,j}$, we have
$$
\|u^e - u^p\|_{L_1(\Omega)}   \le C_1 |\Omega| h^{p+1} \tau  N^{p+1}_{n},
$$
$$
\|u^e - u^p\|_{L_2(\Omega)}  \le C_2 \sqrt{|\Omega|}  h^{p+1} \tau  N^{p+1}_{n},
$$
and, for the cell boundary terms,
$$
(\Sigma_{j=0}^{m-1} |u^e(x_{j+1/2}^-)-u^p(x_{j+1/2}^-)|^2 )^{1/2} \le C_3 \sqrt{|\Omega|} h^{p+1/2}  \tau N^{p+1}_{n}.
$$
\end{lemma}

Next we look into the difference $\tilde{u} - u^e$. In the cell
$\Omega_j$, at time $t_n + \tau$, both of these entropy solutions
$\tilde{u}$ and $u^e$ depend on their initial value in $\Omega_{j-1}
\cup \Omega_j$. More precisely, since $\beta = \max f'(u)$ is the
maximum of wave speed, both entropy solutions $\tilde{u}$ and $u^e$
depend on their initial value in $[x_{j-1/2}-\beta\tau,x_{j-1/2}]
\cup \Omega_j$. Notice that the initial value of these two solutions
are the same in $\Omega_j$ and the difference of their initial
values in $[x_{j-1/2}-\beta \tau,x_{j-1/2}]$ is
$$
 u^e_{j}(t_n)  - u^c_n = J_{n,j}^0 + J_{n,j}^1 (x-x_{j-1/2}) + \cdots +    J_{n,j}^p
 (x-x_{j-1/2})^p/{p!} - R_{n,j}(x).
$$
If $D_{n,j}^0, D_{n,j}^1, \cdots, D_{n,j}^p$ are bounded by
$\tilde{D}$, $\tau \le \gamma  h^{1+\alpha}$, $j>0$, and $x \in
[x_{j-1/2}-\beta\tau,x_{j-1/2}]$, then
\begin{eqnarray}
| u^e_{j}(t_n) - u^c_n| &\le& |D_{n,j}^0 h^{p+1+\mu} + D_{n,j}^1
h^{p+1+\mu-1-\alpha} \beta \gamma h^{1+\alpha} + \cdots
 + D_{n,j}^p h^{p+1+\mu-p(1+\alpha)} (\beta \gamma)^p h^{(1+\alpha)p} /{p!} \, | \nonumber \\
&\le& \tilde{D} h^{p+1+\mu} [1 + \beta \gamma + \cdots + (\beta
\gamma)^p/{p!}] \nonumber\\
&\le&  \tilde{D} e^{\beta \gamma} h^{p+1+\mu} .  \label{disError}
\end{eqnarray}
For $j=0$, there is the the extra residual term $R_{n,0}(x)$, so we
have
\begin{eqnarray}
| u^e_{0}(t_n) - u^c_n| &\le& \tilde{D} h^{p+1+\mu} [1 + \beta
\gamma +
\cdots + (\beta\gamma)^p/{p!}] + \bar{D} (\beta \gamma h^{1+\alpha})^{p+1} / (p+1)! \nonumber \\
&\le&  \tilde{D} e^{\beta \gamma} h^{p+1+\mu},  \label{disError0}
\end{eqnarray}
if $\bar{D} h^{(p+1)\alpha-\mu} \le \tilde{D}$ (which is easy to
satisfy).

 Now, by Theorem 16.1 in the textbook \cite{smoller} by Joel
Smoller, we have
\begin{lemma} \label{lemma3}
If $\beta \tau \le h$ and  $\tau \le \gamma  h^{1+\alpha}$, then
$$
\|\tilde{u}(t_n+\tau) - u^e(t_n+\tau)\|_{L_1(\Omega_j)} \le
\| u^c_n - u^e_{j}(t_n) \|_{L_1 [x_{j-1/2}-\beta\tau,x_{j-1/2}]} \le
(\beta \tau) (\tilde{D} e^{\beta \gamma} h^{p+1+\mu}) .
$$
Consequently,
$$
\|\tilde{u}(t_n+\tau) - u^e(t_n+\tau)\|_{L_1(\Omega)} \le \tau
h^{p+\mu} \beta \tilde{D} e^{\beta \gamma} |\Omega|.
$$
\end{lemma}

{\bf Remarks.} Again, a few remarks may help. (1) Lemma \ref{lemma3}
takes the inter-cell technical discontinuities of the numerical
solution  into account. Obviously, $\tilde{D}$ is playing the role
of a smoothness measurement. (2) Under the strengthened CFL
condition, we are able to allow the value of $J^l_{n,j}$ to be of
order $h^{p+1+\mu-l(1+\alpha)}$. As it is shown in the numerical
examples, when the order of the derivative goes up by one, the power
($h^?$) of the jumps goes down by more than one. The strengthened
CFL condition seems to help here in the error control, at least in
the analysis, even if the jumps of the high order derivatives grow
quickly. When the smoothness deteriorates near the formation of a
shock, the strengthened CFL condition may play a role of suppressing
Runge phenomena, to some extent. Such Runge phenomena and their
transport to the downstream should be what causes numerical oscillations. \\

Now we are ready to state and prove the last lemma to estimate $u^p
- u^h$ , then we will conclude this subsection with the  main
theorem to estimate $\tilde{u}(t_{n+1}) - u^h(t_{n+1})$.

\begin{lemma} \label{lemma4}
There is a computable constant $Q_n$, depending on $S^p_n$, such
that
\begin{equation}  \label{projectionError}
\|u^p(t_n+\tau) - u^h(t_n+\tau)\|_{L_1(\Omega)} \le \tau h^{p+\mu}
Q_n.
\end{equation}
\end{lemma}
\noindent {\bf Proof.} In the cell $\Omega_j$, $u^h$ satisfies the semi-discrete DG scheme (\ref{semi}), that is,
\begin{equation}
(u^h_{t}, v)_{\Omega_j} = (f(u^h), v_x)_{\Omega_j}
                        + f(u^h(x_{j-1/2}^-)) v(x_{j-1/2}^+)
                        - f(u^h(x_{j+1/2}^-)) v(x_{j+1/2}^-). \label{semiRepeat}
\end{equation}
Consider the piecewise strong solution $u^e$, which is the
restriction of $u^{e}_j$  in $\Omega_j$.  Multiplying $u^e_t +
f(u^e)_x = 0$ by a test function $v$, integrating in $\Omega_j$,
after using integration by parts, we get
\begin{equation}
(u^e_{t}, v)_{\Omega_j} = (f(u^e), v_x)_{\Omega_j}
                        + f(u^e(x_{j-1/2}^+)) v(x_{j-1/2}^+)
                        - f(u^e(x_{j+1/2}^-)) v(x_{j+1/2}^-). \label{piecewiseSolu}
\end{equation}
Since $(u^p,v)_{\Omega_j} = (u^e,v)_{\Omega_j} $ for all values of $\tau \in [0,\tau_n]$,  $(u^p_t,v)_{\Omega_j} = (u^e_t,v)_{\Omega_j} $. By adding and subtracting terms in (\ref{piecewiseSolu}), we get
\begin{eqnarray*}
(u^p_{t}, v)_{\Omega_j} &=& (f(u^p), v_x)_{\Omega_j}
                        + f(u^e(x_{j-1/2}^+)) v(x_{j-1/2}^+)
                        - f(u^p(x_{j+1/2}^-)) v(x_{j+1/2}^-)    \\
&+&       (f(u^e)-f(u^p), v_x)_{\Omega_j} -  [ f(u^e(x_{j+1/2}^-)) - f(u^p(x_{j+1/2}^-)) ]  v(x_{j+1/2}^-) .
\end{eqnarray*}
Now let $\xi = u^p - u^h$, and let $v = \xi$.
Subtracting the last equation by (\ref{semiRepeat}),  we get
\begin{eqnarray}
(\xi_{t},\xi)_{\Omega_j} &=& (f(u^p)-f(u^h), \xi_x)_{\Omega_j} \nonumber \\
&-& [ f(u^p(x_{j+1/2}^-)) -  f(u^h(x_{j+1/2}^-))] \, \xi(x_{j+1/2}^-)   \nonumber  \\
&+&     [  f(u^e(x_{j-1/2}^+)) -   f(u^h(x_{j-1/2}^-)) ] \, \xi(x_{j-1/2}^+)  \nonumber \\
&+&  (f(u^e)-f(u^p), \xi_x)_{\Omega_j}    \nonumber \\
& -&  [ f(u^e(x_{j+1/2}^-)) - f(u^p(x_{j+1/2}^-)) ] \,
\xi(x_{j+1/2}^-) . \label{xi}
\end{eqnarray}
First we focus on the third line of the last equation.
\begin{eqnarray}
| f(u^e(x_{j-1/2}^+)) -   f(u^h(x_{j-1/2}^-))| &\le& \beta \, | u^e(x_{j-1/2}^+) -   u^h(x_{j-1/2}^-)|   \nonumber \\
&\le&  \beta \, | u^e(x_{j-1/2}^+) -  u^e(x_{j-1/2}^-)| \nonumber \\
&+&   \beta \, |u^e(x_{j-1/2}^-) - u^p(x_{j-1/2}^-)|  \nonumber\\
&+&   \beta \, |u^p(x_{j-1/2}^-) - u^h(x_{j-1/2}^-)|
\label{secondTerm}
\end{eqnarray}
for $j>0$. As for $j=0$, by verifying $u^e(x_{-1/2}^-) =
u^h(x_{-1/2}^-)$ from the boundary conditions, we have
\begin{eqnarray}
| f(u^e(x_{-1/2}^+)) -   f(u^h(x_{-1/2}^-))| &\le& \beta \, | u^e(x_{-1/2}^+) -   u^h(x_{-1/2}^-)|   \nonumber \\
&=&  \beta \, | u^e(x_{-1/2}^+) -  u^e(x_{-1/2}^-)|
\label{secondTerm0} .
\end{eqnarray}
Here, we use a brief notation $ A = u^e(x_{j-1/2}^+)
$, and $ B = u^e(x_{j-1/2}^-)$. Recall that $ u^e = u^e_j$  in
$\Omega_j$ and $ u^e = u^e_{j-1}$  in $\Omega_{j-1}$. Also recall
that  $ u^e_j$ is extended to $\R_{n,j}$. By the characteristic
line theory, we know that
$$
A = u^e_j(t_n+\tau, x_{j-1/2}) = u^e_j(t_n,x_{j-1/2} - \tau f'( u^e_j(t_n+\tau, x_{j-1/2})) )
= u^e_{j}(t_n, x_{j-1/2} - \tau f'(A))
$$
and
$$
B = u^e_{j-1}(t_n+\tau, x_{j-1/2}) = u^e_{j-1}(t_n,x_{j-1/2} - \tau f'( u^e_{j-1}(t_n+\tau, x_{j-1/2})) )
= u^e_{j-1}(t_n,x_{j-1/2} - \tau f'(B))  .
$$
So
\begin{eqnarray*}
|A-B| &=&    |u^e_{j}(t_n,x_{j-1/2} - \tau f'(A))  -  u^e_{j-1}(t_n,x_{j-1/2} - \tau f'(B))  | \\
         &\le&  |u^e_{j}(t_n,x_{j-1/2} - \tau f'(A))  -  u^e_{j}(t_n,x_{j-1/2} - \tau f'(B))  | \\
           &+&  |u^e_{j}(t_n,x_{j-1/2} - \tau f'(B))  -  u^e_{j-1}(t_n,x_{j-1/2} - \tau f'(B))  |
\end{eqnarray*}
By Lemma \ref{lemma1}, the first spatial derivative of
$|u^e_{j}(t_n)|$ is bounded by $N^1_{n,j}$. Assume that $|f''|$ is
bounded by $\delta<\infty$. Since $f'(B) \le \beta$, for
sufficiently small $\tau$, according to the inequalities
(\ref{disError}) and (\ref{disError0}),
\begin{eqnarray*}
|u^e_{j}(t_n,x_{j-1/2} - \tau f'(B)) - u^e_{j-1}(t_n,x_{j-1/2} - \tau f'(B))| \le \tilde{D} e^{\beta \gamma} h^{p+1+\mu}.
\end{eqnarray*}
Therefore,
$$
|A-B| \le N^1_{n,j} \delta \tau |A-B| +  \tilde{D} e^{\beta \gamma} h^{p+1+\mu}.
$$
Due to the fact that $\tau$ is very small, we have
\begin{equation} \label{ee}
| u^e(x_{j-1/2}^+) -  u^e(x_{j-1/2}^-)| = |A-B| \le (1+2N^1_{n,j}
\delta \tau)  \tilde{D} e^{\beta \gamma} h^{p+1+\mu}.
\end{equation}
Now, plug (\ref{secondTerm}) into (\ref{xi}), take the sum over all cells, we get
\begin{eqnarray}
(\xi_{t},\xi)_{\Omega} =&& (f(u^p)-f(u^h), \xi_x)_{\Omega} \nonumber \\
&-& \Sigma_{j=0}^{m-1} [ f(u^p(x_{j+1/2}^-)) -  f(u^h(x_{j+1/2}^-))] \, \xi(x_{j+1/2}^-)   \nonumber  \\
&+&    \Sigma_{j=0}^{m-1} [  f(u^e(x_{j-1/2}^+)) -   f(u^h(x_{j-1/2}^-)) ] \, \xi(x_{j-1/2}^+)  \nonumber \\
&+&   (f(u^e)-f(u^p), \xi_x)_{\Omega} \nonumber  \\
& -&  \Sigma_{j=0}^{m-1} [ f(u^e(x_{j+1/2}^-)) - f(u^p(x_{j+1/2}^-)) ] \, \xi(x_{j+1/2}^-)  \nonumber \\
\le&& \beta \, \| \xi\|_{L_2(\Omega)} \, \| \xi_x \|_{L_2(\Omega)} \nonumber \\
&+&  \beta \, \Sigma_{j=0}^{m-1} \xi^2(x_{j+1/2}^-) \nonumber\\
&+& \beta \, \Sigma_{j=0}^{m-1} | u^e(x_{j-1/2}^+) -  u^e(x_{j-1/2}^-)| \, | \xi(x_{j-1/2}^+) | \nonumber \\
&+& \beta \, \Sigma_{j=1}^{m-1} |u^e(x_{j-1/2}^-) - u^p(x_{j-1/2}^-)| \, |\xi(x_{j-1/2}^+)| \nonumber \\
&+& \beta \, \Sigma_{j=1}^{m-1} |u^p(x_{j-1/2}^-) - u^h(x_{j-1/2}^-)| \, |\xi(x_{j-1/2}^+)| \nonumber \\
&+&  \beta \, \|u^e-u^p\|_{L_2(\Omega)} \, \| \xi_x \|_{L_2(\Omega)} \nonumber \\
&+&  \beta \, \Sigma_{j=0}^{m-1} |u^e(x_{j+1/2}^-) -
u^p(x_{j+1/2}^-)| \, |\xi(x_{j+1/2}^-)| .\nonumber
\end{eqnarray}
By using (\ref{ee}), the estimates on the projection error given in
Lemma \ref{lemma2}, and the standard inverse inequalities
(\cite{qz2009}, section 3.3), we can get computable constants $C_4,
C_5$ and $C_6$, such that
$$
\frac{d}{dt} \|\xi \|_{L_2(\Omega)} \le \frac{C_4}{h} \|\xi
\|_{L_2(\Omega)} + \frac{C_5}{h}  h^{p+1+\mu} + \frac{C_6}{h} \tau
h^{p+1} .
$$
Integrating the last differential inequality, noticing the fact that
$\xi = 0$ at $\tau = 0$, also noticing that $\tau \le \gamma h^{1+\alpha} \le \gamma h^\mu$, we have a computable constant $C_7$, such
that
$$
 \|\xi \|_{L_2(\Omega)} \le C_7 \tau h^{p+\mu}.
$$
Since $\| \xi \|_{L_1(\Omega)} \le \sqrt{|\Omega|} \| \xi
\|_{L_2(\Omega)}$, we have (\ref{projectionError}) and Lemma
\ref{lemma4} proven.
\# \\

Combining Lemma \ref{lemma3}, Lemma \ref{lemma2} and Lemma \ref{lemma4}, we have the following theorem.
\begin{theorem} \label{spaceTheorem}
There is a computable constant $\F(S^p_n)$, depending on the flux
function $f$, the known constants of the interpolation/projection
error estimates, the known constants of the inverse inequalities,
and the components of the spatial smoothness indicator $S^p_n$, such
that
$$
\|\tilde{u}(t_{n+1})  - u^h(t_{n+1})\|_{L_1(\Omega)} \le \tau_n
h^{p+\mu} \F(S^p_n).
$$
\end{theorem}

\subsection{Estimating $u^h(t_{n+1}) - u^c_{n+1}$, proof of
(\ref{timeError}) }

The temporal smoothness indicator $T^k_n$ informs us about the
boundedness of the temporal derivatives of $u^h$
at $t=t_n$. We need to make sure that the boundedness of $T^k_n$
can guarantee the boundedness of the temporal derivatives of
$u^h$ for all $t \in [t_n,t_{n+1}]$.

\begin{lemma} \label{lemma5}
There is a computable constant K, depending on the spatial
smoothness indicator $S^p_n$, such that
\begin{equation} \label{uh}
\| u^h(t_n+\tau) \| _{L_{\infty}(\Omega)} \le \| u^c_n \|
_{L_{\infty}(\Omega)} + K h.
\end{equation}
For each integer $l \in \{1,\cdots, k+1\}$, there is a pair of
computable constants $c_l$ and $d_l$, such that, for all $t \in
[t_n,t_{n+1}]$,
\begin{equation} \label{uhl}
\| \frac{\partial^l}{\partial t ^l} u^h \|
_{L_{\infty}(\Omega)} \le (1+c_l h^{\alpha}) \|
\frac{\partial^l}{\partial t ^l} u^h( t_n )\|
_{L_{\infty}(\Omega)} + d_l h^{\alpha}.
\end{equation}
For each $l$, $c_l$ and $d_l$ only depend on the $L_{\infty}$-norms
of the lower order derivatives.
\end{lemma}
\noindent {\bf Proof.}
\begin{eqnarray*}
\|  u^h(t_n+\tau) \| _{L_{\infty}(\Omega)}
&\le& \|\tilde{u}(t_n+\tau) \| _{L_{\infty}(\Omega)} \\
&+& \| \tilde{u}(t_n+\tau) - u^e(t_n+\tau)\|
_{L_{\infty}(\Omega)} \\
&+& \| u^e(t_n+\tau) - u^p(t_n+\tau)\|
_{L_{\infty}(\Omega)} \\
&+& \| u^p(t_n+\tau) - u^h(t_n+\tau)\|
_{L_{\infty}(\Omega)} .
\end{eqnarray*}
$\|\tilde{u} \| _{L_{\infty}(\Omega)}$ is bounded by
$\|u^c_n\| _{L_{\infty}(\Omega)}$, because the maximum of the
entropy solution $\tilde{u}$ does not increase. The
smallness of $\| u^p- u^h\|
_{L_{\infty}(\Omega)}$ can be obtained from Lemma \ref{lemma4} and
an application of the inverse inequality. The smallness of $\|
u^e - u^p\|
_{L_{\infty}(\Omega)}$ is given by Lemma \ref{lemma2}. The smallness
of $\| \tilde{u} - u^e\|
_{L_{\infty}(\Omega)}$ can be obtained by the same method used in
proving (\ref{ee}). Consequently, we have the estimate (\ref{uh}).

In the proof of (\ref{uhl}), let's set the notation $z = z(\tau) =
u^h(t_n+\tau)$, and $z^{(l)} = \frac{\partial^l}{\partial t
^l} z$. By differentiating the semi-discrete DG scheme
$$
(z_{t}, v)_{\Omega_j} = (f(z), v_x)_{\Omega_j}
                        + f(z(x_{j-1/2}^-)) v(x_{j-1/2}^+)
                        - f(z(x_{j+1/2}^-)) v(x_{j+1/2}^-)
$$
with respect to $t$, we get
\begin{eqnarray*}
(z^{(1)}_{t}, v)_{\Omega_j} &=& (f'(z) z^{(1)}, v_x)_{\Omega_j} \\
        &+& f'(z(x_{j-1/2}^-)) z^{(1)}(x_{j-1/2}^-) v(x_{j-1/2}^+) \\
        &-& f'(z(x_{j+1/2}^-)) z^{(1)}(x_{j+1/2}^-) v(x_{j+1/2}^-),
\end{eqnarray*}
\begin{eqnarray*}
(z^{(2)}_{t}, v)_{\Omega_j} &=& (f'(z) z^{(2)}, v_x)_{\Omega_j}
        + (f''(z) [z^{(1)}]^2 , v_x)_{\Omega_j} \\
        &+& f'(z(x_{j-1/2}^-)) z^{(2)}(x_{j-1/2}^-) v(x_{j-1/2}^+)
        + f''(z(x_{j-1/2}^-)) [z^{(1)}(x_{j-1/2}^-)]^2 v(x_{j-1/2}^+)\\
        &-& f'(z(x_{j+1/2}^-)) z^{(2)}(x_{j+1/2}^-) v(x_{j+1/2}^-)
        - f''(z(x_{j+1/2}^-)) [z^{(1)}(x_{j+1/2}^-)]^2 v(x_{j+1/2}^-),
\end{eqnarray*}
and similar equations for $z^{(l)}$, $l=3, \cdots, k+1$. It is easy
to observe that the equation for $z^{(l)}$ is linear on $z^{(l)}$,
$l=1,2,\cdots,k+1$. Moreover, it depends on the derivatives of the
flux function $f$ and and products of lower order derivatives of
$z$.

In order to estimate the $L_{\infty}$-norm of $z^{(1)}$, we expand
$z^{(1)}$ by the normalized Legendre polynomial basis functions $\{
\phi_{j,i}: i=0,\cdots,p \}$ of $V_h$ in each cell $\Omega_j$.
$(\phi_{j,i}, \phi_{j,i})_{\Omega_j} = h/2$.
$$
z^{(1)}(\tau)|_{\Omega_j} = \Sigma_{i=0}^p q^{(1)}_{j,i}(\tau) \phi_{j,i}.
$$
Under this basis, let $q^{(1)}$ be the vector consisting of all the
$q^{(1)}_{j,i}$, $j=0,\cdots,m-1$; $i=0,\cdots,p$. We can rewrite
the equation for $z^{(1)}$ as
$$
\frac{h}{2} \frac{d}{dt} q^{(1)} = A^{(1)}(\tau) q^{(1)},
$$
where $A^{(1)}(\tau)$ is the matrix obtained from the righthand side
of the equation about $z^{(1)}$. The entries of $A^{(1)}(\tau)$
depend on the wave speed $f'(z)$. Since $\|  z \|
_{L_{\infty}(\Omega)}$ is bounded according to (\ref{uh}) and $f$ is
smooth, the entries of $A^{(1)}(\tau)$ are bounded. Besides, the
entries does not depend on $h$, and there is at most $2p+2$ entries
in each row of $A^{(1)}(\tau)$ not equal to zero. Solving for
$q^{(1)}$ from the last ODE, we get
$$
q^{(1)}(\tau) = e^{\frac{2}{h} \int_0^{\tau} A^{(1)}(\tau) d\tau}
q^{(1)}(0).
$$
From this solution, it is easy to see that, there is a constant
$\tilde{A}$ depending on the entries of $A^{(1)}(\tau)$, such that
$$
\|q^{(1)}(\tau) - q^{(1)}(0)\|_{\infty} \le \frac{\tilde{A}\tau}{h}
\|q^{(1)}(0)\|_{\infty} \le \tilde{A} \gamma h^{\alpha}
\|q^{(1)}(0)\|_{\infty} .
$$
Due to the equivalence of $\|z^{(1)}\|_{L_{\infty}(\Omega)}$ and
$\|q^{(1)}\|_{\infty}$, we have constants $\tilde{B} =
(p+1)\sqrt{(2p+1)/2}$ and $\tilde{C} = \sqrt{2}$, such that
\begin{eqnarray*}
\|z^{(1)}(\tau)\|_{L_{\infty}(\Omega)} &\le&
\|z^{(1)}(0)\|_{L_{\infty}(\Omega)} +
\|z^{(1)}(\tau) - z^{(1)}(0) \|_{L_{\infty}(\Omega)} \\
&\le& \|z^{(1)}(0)\|_{L_{\infty}(\Omega)} +
\tilde{B} \|q^{(1)}(\tau) - q^{(1)}(0) \|_{\infty} \\
&\le& \|z^{(1)}(0)\|_{L_{\infty}(\Omega)} +
\tilde{B} \tilde{A} \gamma h^{\alpha} \|q^{(1)}(0) \|_{\infty} \\
&\le& \|z^{(1)}(0)\|_{L_{\infty}(\Omega)} + \tilde{C} \tilde{B}
\tilde{A} \gamma h^{\alpha} \|z^{(1)}(0) \|_{L_{\infty}(\Omega)}
\end{eqnarray*}
This proves (\ref{uhl}) for $l=1$, with $c_1 = \tilde{C}\tilde{B}
\tilde{A} \gamma$ and $d_1=0$. For $l \ge 2$, one can carry out the
proof in
the same way.  \# \\

{\bf Remarks.} (1) Lemma \ref{lemma5} confirms that one can
essentially use the value of $\frac{\partial^l}{\partial t ^l}
u^h(0, t_n, u^c_n)$ as a estimate of $\frac{\partial^l}{\partial t
^l} u^h(\tau, t_n, u^c_n)$ in the entire time step $[t_n,t_{n+1}]$.
(2) The result of Lemma \ref{lemma5} serves our purpose of local
smoothness validation. However, neither the result nor the method of
proof can/should be generalized to long term, because no numerical
diffusion is taken into account.\\

Based on the boundedness of the temporal derivatives proven in Lemma
\ref{lemma5}, it is trivial to conclude with the next theorem.

\begin{theorem}
There is a computable function $\G(T^k_n,S_n^{p})$, such that
\begin{equation}
\|u^h(t_{n+1}) - u^c_{n+1}\|_{L_1(\Omega)}
\le \tau_n^{k+1} \G(T^k_n,S_n^{p}) .
\end{equation}
\end{theorem}

\section{Numerical evidences}

From the error estimation inequalities (\ref{spaceError}) to
(\ref{global}), once we show the boundedness of all the components
of the smoothness indicators, the rest of the error estimates is
essentially {\it a priori}. Therefore, in order to demonstrate that
our analysis works, it suffices to display the computed smoothness
indicators.\\

\begin{figure}
  \includegraphics[width=6.0in,height=2.2in]{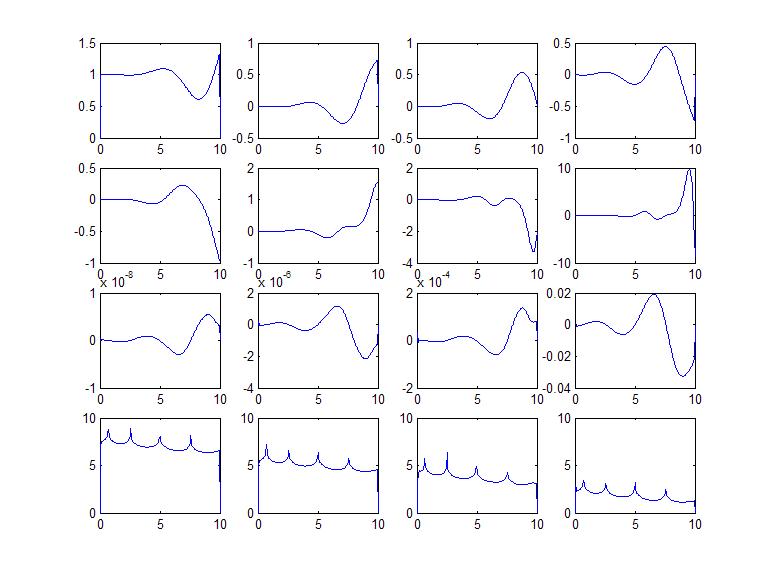}
  \caption{Smoothness Indicators, t=0.05}\label{si010}
\end{figure}

\begin{figure}
  \includegraphics[width=6.0in,height=2.2in]{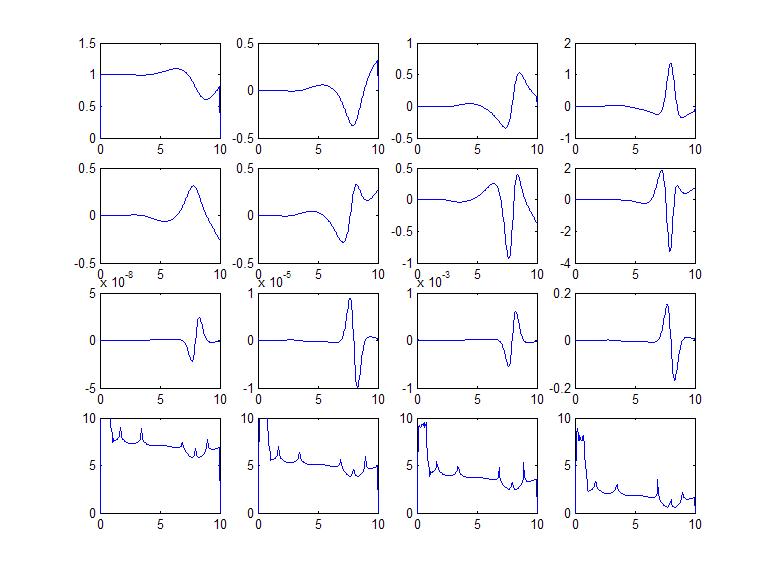}
  \caption{Smoothness Indicators, t=1.05}\label{si210}
\end{figure}

\begin{figure}
  \includegraphics[width=6.0in,height=2.2in]{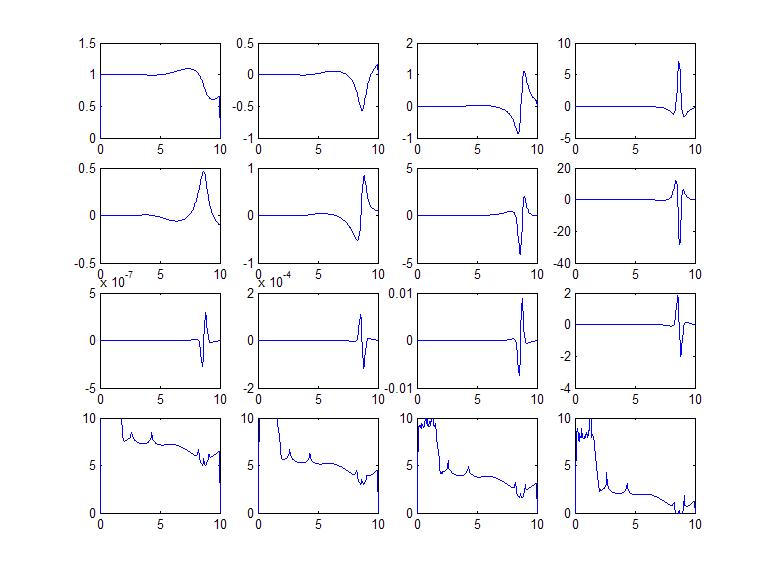}
  \caption{Smoothness Indicators, t=2.0}\label{si395}
\end{figure}

\noindent {\bf Example 1.} In the first example, we solve  Burgers'
equation
$$
u_t + (u^2/2)_x = 0
$$
with the boundary condition $u_L(t) = 1$ and initial condition
$$
u_I(x) = 1-(x/11)^3 \sin(x)
$$
in $x \in \Omega = [0,10]$. In this numerical example with $p=3$ and
$k=3$, the cell size is $h=0.05$, while the time step size is
$\tau_n = 0.005$.. The solution has been computed in $t \in [0,2]$.
The smoothness indicators at $t=0.05$, $t=1.05$, and $t=2.0$ are
shown in Figure \ref{si010}, Figure \ref{si210},  and Figure
\ref{si395} respectively.

In each figure, the four plots in the top row are $M^0_n (=u^c_n)$,
$M^1_n$, $M^2_n$, and $M^3_n$, from left to right. The index $j$ is
dropped because each curve contains the values of $M^l_{n,j}$ for
$j=0,\cdots, 199$. The four plots in the second row are the temporal
smoothness indicators $u^h_t(t_n)$, $u^h_{tt}(t_n)$, $u^h_{ttt}(t_n)$, and $u^h_{tttt}(t_n)$.
The four plots in the third row are the jumps $J^0_n$, $J^1_n$,
$J^2_n$, and $J^3_n$. In order to view the jumps from a better
perspective, we show $\log_h |J^0_n|$, $\log_h |J^1_n|$, $\log_h
|J^2_n|$ and $\log_h |J^3_n|$ in the fourth row. Since $J^l_{n,j} =
D^l_{n,j} h^{p+1+\mu-l(1+\alpha)}$, the plot of $\log_h |J^l_{n,j}|
= p+1+\mu-l(1+\alpha) + \log_h |D^l_{n,j}|$ reveals the order
($h^?$) of the jumps. Since $p, \mu$ and $\alpha$ are all known, the
values of $D^l_{n,j}$ can be computed. Consequently, we can find
$\tilde{D}$.


It is easy to see the boundedness of $M^l_n$ and $T^3_n$ in the
figures when/where the solution is smooth. It is also easy to see
that the order of the jumps $J^l_n$ is as expected in the error
analysis, or even smaller. These observations are sufficient to
support the error estimates given in the paper.

In addition, we have also observed some interesting phenomena. (1)
$log_h|J^0_n| - log_h|J^1_n| \approx 2$, $log_h|J^1_n| -
log_h|J^2_n| \approx 1.4$, $log_h|J^2_n| - log_h|J^3_n| \approx
1.8$. There seems to be something related to the odd or even degrees
of the polynomials. (2) Long before the formation of a shock
($t=2.0$, $u_x \ge -0.6$), the fourth and third derivatives have
grown significantly in a very narrow subdomain. The approximation
benefit of the higher degree polynomials and the high order
Runge-Kutta scheme will soon be lost locally at the spot. It seems
that adaptive treatments need to kick in early. If not, there will
be ``numerical instability" showing up, ruining the numerical
solution.

\begin{figure}
  \includegraphics[width=6.0in,height=2.2in]{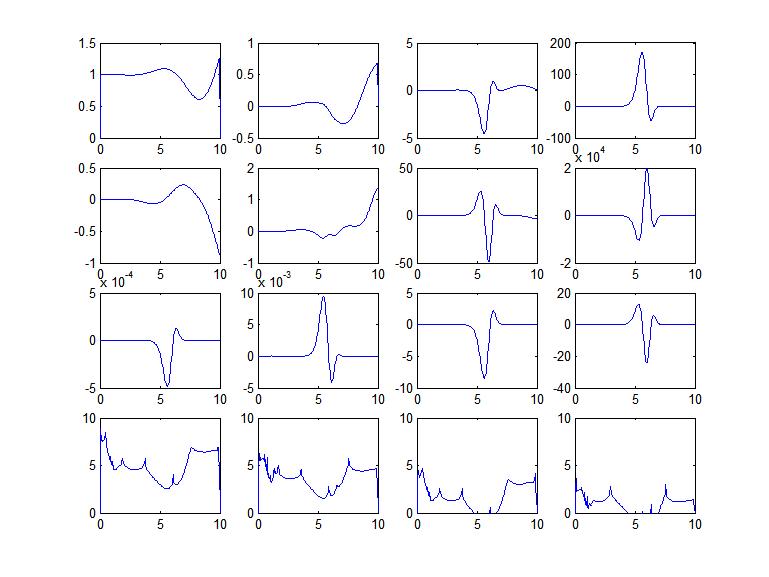}
  \caption{Smoothness Indicators, $\tau_n = 0.0075, t=0.12$}\label{sibad}
\end{figure}

In Figure \ref{sibad}, we show another numerical solution of the
same problem, computed with $h=0.05$ (same as before) and $\tau_n =
0.0075$ (50\% larger). The plots are made at $t=0.12$, after 16 time
steps from the initial time $t=0$. With the improperly increased
time step size, although the solution (presented by $M^0_n$ in the
upper left corner plot) itself has not obviously shown anything
wrong from the point of view of numerical stability (boundedness of
solution, TVD, etc.), the higher order derivatives and jumps in the
indicators have been increased significantly. The explanation is
that the RKDG scheme for this problem with $(p,k,h,\tau_n) =
(3,3,0.05, 0.0075)$ does not maintain numerical smoothness. As a
consequence, the optimal approximation order must have been lost.
The example seems to indicate the following: the strengthened CFL
condition and the numerical diffusion from the Godunov flux are
needed not only for numerical stability, but also for numerical
smoothness maintenance. More attention should be paid to numerical
smoothness when we are concerned with high order error estimates.

The smoothness indicators can be used to diagnose the loss of
numerical smoothness in an early stage, before too much damage is
done to the global error. Of course, an algorithm needs to be
designed for such diagnoses. We did run a separate case: after the
first 5 steps at $\tau_n = 0.0075$, $\tau_n$ is reduced back to
$0.005$. The spurious mode created in the first 5 steps were
repaired in the following steps of smaller size. Nevertheless, the
damage to the global error is done, unless we redo it. Further
investigation in this direction can help in finding an optimal time
step size. \\

\begin{figure}
  \includegraphics[width=6.0in,height=2.5in]{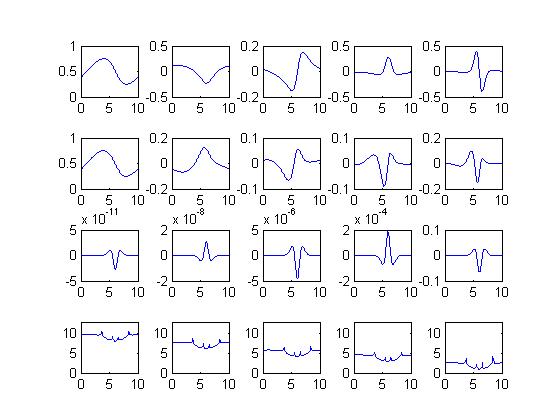}
  \caption{Smoothness Indicators, $p=4, \tau_n = 0.005, t=1.0$}\label{p200}
\end{figure}

\noindent {\bf Example 2.} In the second example, we show the
solution of the Burgers' equation on $[0,10]$ with the initial
condition
$$
u_I(x,0) = \frac{1}{2} + \frac{1}{4}\sin(\pi x /5)
$$
and the periodic boundary condition. $k=3$, $p=4$, $h=0.05$,
$\tau_n=0.005$. Figure \ref{p200} shows the numerical solution and
its smoothness indicators at $t=1$, when it is still far from any
shock formation.

In Figure \ref{p200}, the five plots in the top row are $M^0_n
(=u^c_n)$, $M^1_n$, $M^2_n$, $M^3_n$, and $M^4_n$, from left to
right. The five plots in the second row are the temporal smoothness
indicators $u^c_n$, $u^h_t(t_n,)$, $u^h_{tt}(t_n)$,
$u^h_{ttt}(t_n)$, and $u^h_{tttt}(t_n)$. The five
plots in the third row are the jumps $J^0_n$, $J^1_n$, $J^2_n$,
$J^3_n$ and $J^4_n$. In the fourth row, we have $\log_h |J^0_n|$,
$\log_h |J^1_n|$, $\log_h |J^2_n|$, $\log_h |J^3_n|$ and $\log_h
|J^4_n|$. The values of the indicators are again what we expected
and what we need to support the analysis. Obviously, the scheme has
maintained the smoothness of the numerical solution, which
guarantees that the local error of the next time step will be of
optimal order.

\section{Conclusion remarks}

\noindent {\bf A. Choice of norm for error propagation analysis}

We prefer to use the $L_1$-norm for error propagation analysis
because of the well-known $L_1$-contraction property. Other than the
$L_1$-contraction, a typical error propagation rate estimate for a
time step contains a growth factor of the form $1+C\tau$.  If we
choose the $L_2$-norm for error propagation, it is easy to show that
the constant $C$ is proportional to $\sqrt{|u_x f''(u)|}$. If we use
numerical error propagation instead of PDE's error propagation, $C$
will become bigger.  ``Bigger by how much" depends on the complexity
of a numerical scheme.  The appearance of $u_x f''(u)$ in the
$L_2$-norm error propagation rate estimate implies that $L_2$-norm
error propagation analysis based on ``worst case scenario" cannot be
generalized to solutions with a shock or near a shock. Since large
local error is expected to appear around the self-sharpening of a
solution, the real scenario of a numerical solution is probably very
close to the ``worst case scenario". $L_1$-norm error propagation analysis does not have this difficulty. \\

\noindent {\bf B. How to deal with shocks and contact
discontinuities?}

When there is a shock or contact discontinuity, it will be detected
by the smoothness indicators, as shown in \cite{rumsey}. Certain
quantitative criteria need to be developed to determine what kind of
discontinuity is present according to the behavior of the
indicators. It is also needed to determine if the discontinuity is
well-caught, or some level of numerical ``instability" has occurred.
A decision should be made on the treatment of the discontinuity,
including the use of a limiter or a local front tracking technique.
After all of these have been done, we can consider error estimation.
Error propagation is still to be estimated by using
$L_1$-contraction. Within each time step, in the smooth pieces of
the solution, we can apply the error estimates given in this paper.
At the discontinuities, we have to estimate the error according to
the scheme. It is nice that the complexity of local error analysis
does not get into the error propagation of the PDE.\\

\noindent {\bf C. The process of sharpening before shock formation
may be most difficult}

It might be the hardest to estimate error where a shock is forming
but not yet fully developed. In this relatively wide space-time
region, the solution's high order derivatives have become larger,
causing difficulties for approximation. Adaptive algorithms need to
be designed, and employed according to the smoothness indicators.
As seen in Figure \ref{si395} of the first numerical example, the
smoothness indicators can find the local sharp growth of the higher
order derivatives and their jumps. The logarithm plots of the jumps have
shown a clear exclusive pattern for a point of future shock. \\

\noindent {\bf D. Generalization to multi-dimensional problems}

We checked the proofs to the end of generalizing the results to 2-D
scalar conservation laws. It seems to us that such a generalization
should not meet any major difficulty. Generalization to hyperbolic
systems will face the lack of $L_1$-contraction. \\

\noindent {\bf E. {\it a posteriori} vs. {\it a priori} estimates}

The error analysis of this work is {\it a posteriori} because we
depend on the computed smoothness indicators to compute the error
estimates. However, if one can prove the boundedness of these
smoothness indicators in advance, the error estimates can be
converted to {\it a priori} error bounds. In this sense, under the
concept of numerical smoothness, {\it a priori} and {\it a
posteriori} error analysis has been united in the same framework.
Moreover, our estimates are {\it a posteriori} in the sense that the
smoothness indicators $S^p_n$ and $T^k_n$ are computed after $u^c_n$
has been obtained. As for the time step $[t_n,t_{n+1}]$, the
smoothness indicators needed for the local error estimates of the
step are computed before the local computation toward $u^c_{n+1}$
has started. In this sense, our error estimation is locally {\it a
priori}, which will be more efficient
if adaptive treatments are desired.  \\

\noindent {\bf F. Numerical smoothness of RKDG}

In the error analysis, we actually depend on the smoothness
indicators to provide the needed numerical diffusion. That is, we
take advantage of the RKDG method to include the needed numerical
smoothness maintenance into the error analysis. The original
designers of the scheme should get the credit for inventing a scheme
with such properties. Since the numerical smoothness indicators
$S^p_n$ and $T^k_n$ are computed at $(t_n, u^c_n)$, Lemma
\ref{lemma1} and Lemma \ref{lemma5} are needed to establish the
smoothness of $\tilde{u}$ and $u^h$ for $t \in [t_n,t_{n+1}]$. Lemma
\ref{lemma1} shows the local smoothness preserving property of the
PDE's strong solutions (in a special case useful for the analysis).
Lemma \ref{lemma5} shows the local smoothness preserving property of
the semi-discrete scheme. We only need these local smoothness proofs
because smoothness is only needed in dealing with local error
estimates.

\bibliographystyle{amsplain}

%\end{article}

\end{document}